%% file: root.tex
\documentclass[letterpaper, 10 pt, conference]{ieeeconf}  

\IEEEoverridecommandlockouts                              

\usepackage{times} 
\usepackage{amsmath} 
\usepackage{amssymb}  
\usepackage{amsfonts}
\usepackage{mathtools}
\usepackage{algorithm}
\usepackage{algpseudocode}
\usepackage{subfiles}
\usepackage{xcolor}
\usepackage{booktabs}
\usepackage{multirow}
\usepackage{subcaption}
\usepackage{placeins}
\usepackage[noadjust]{cite}
\usepackage{hyperref}

\newtheorem{problem}{Problem}

\title{\LARGE \bf
Distributed Stochastic Search for Multi-Agent Model Predictive Control
}

\author{Taehyun Yoon$^{1}$, 
Augustinos D. Saravanos$^{2}$, and 
Evangelos A. Theodorou$^{1}$
\thanks{This work was supported by ARO Award \#W911NF2010151.  
 Taehyun Yoon acknowledges support by the Fulbright Scholarship.  
 Augustinos Saravanos acknowledges support by the A. Onassis Foundation Scholarship.
 }
\thanks{$^{1}$Taehyun Yoon and Evangelos A. Theodorou are with the Daniel Guggenheim School of Aerospace Engineering, Georgia Institute of Technology, Atlanta GA, 30332, USA. {\tt\small \{taehyun.yoon, \,\, evangelos.theodorou\}@gatech.edu}}%
\thanks{$^{2}$Augustinos D. Saravanos was with the School of Electrical and Computer Engineering, Georgia Institute of Technology, Atlanta, GA 30332, USA, during his involvement in this work. He is now with the Department of Aeronautics and Astronautics, Massachusetts Institute of Technology, Cambridge, MA 02139, USA. {\tt\small asaravan@mit.edu}}%
}

\begin{document}
\include{commands}

\maketitle
\thispagestyle{empty}
\pagestyle{empty}

\begin{abstract}
Many real-world multi-agent systems exhibit nonlinear dynamics and complex inter-agent interactions.  
As these systems increase in scale, the main challenges arise from achieving scalability and handling nonconvexity.  
To address these challenges, this paper presents a distributed sampling-based optimization framework for multi-agent model predictive control (MPC).  
We first introduce stochastic search, a generalized sampling-based optimization method, as an effective approach to solving nonconvex MPC problems because of its exploration capabilities.  
Nevertheless, optimizing the multi-agent systems in a centralized fashion is not scalable as the computational complexity grows intractably as the number of agents increases.  
To achieve scalability, we formulate a distributed MPC problem and employ the alternating direction method of multipliers (ADMM) to leverage the distributed approach.  
In multi-robot navigation simulations, the proposed method shows a remarkable capability to navigate through nonconvex environments, outperforming a distributed optimization baseline using the interior point optimizer (IPOPT).
In a 64-agent multi-car formation task with a challenging configuration, our method achieves 100\% task completion with zero collisions, whereas distributed IPOPT fails to find a feasible solution.
\end{abstract}

\section{\textsc{Introduction}}
Model predictive control (MPC) is one of the successful paradigms adopted in various systems, including aerospace control systems \cite{eren2017model}, chemical engineering \cite{kumar2012model}, and autonomous vehicles \cite{garikapati2024autonomous}.  
In multi-agent settings, MPC is an attractive approach that supports reactive communication and coordination between agents under tight real-time budgets.  
The use cases of multi-agent MPC cover transportation networks \cite{de2010multi}, power networks \cite{negenborn2007multi}, multi-robot systems \cite{dai2017distributed}, and many other domains.  

Many real-world systems exhibit highly complex dynamics that cannot be fully captured by linear models.  
To control nonlinear systems, MPC problem can be cast as a sequence of nonlinear programming (NLP) problems for each time step.  
However, NLP problems for controlling nonlinear systems are inherently nonconvex \cite{rawlings2020model}.  
For nonconvex problems, there is no guarantee that the global minimum can be achieved \cite{boyd2004convex}, and convex optimization methods are prone to convergence to an undesired local solution \cite{tenny2004closed}.  

Beyond convex optimizers, nonconvexity in nonlinear control problems is addressed by: i) sampling-based optimization; ii) convexification \cite{schulman2014motion}; and iii) mixed integer programming formulation \cite{richards2002spacecraft}.    
In particular, sampling-based optimization methods \cite{williams2017information, rubinstein1997optimization} have shown promise in solving nonconvex nonlinear MPC problems in real time due to their parallelization and exploration capabilities.  
In this context, sampling-based optimization is an appropriate candidate to solve nonconvex nonlinear MPC problems.  


Distributed optimization architectures based on the alternating direction method of multipliers (ADMM) \cite{boyd2011distributed} have been receiving increasing interest to achieve scalability to large multi-agent systems in control and autonomy \cite{saravanos2023distributed, Saravanos-RSS-21, rey2018fully, ferranti2022distributed}. Naturally, several works have integrated distributed ADMM into multi-agent MPC with significant success.
Specifically, Summers and Lygeros \cite{summers2012distributed} used ADMM to solve the dual decomposition of the distributed MPC (DMPC) problem, where the problem is assumed to be convex.  
Saravanos et al. \cite{saravanos_2024distributedmpcs} proposed an ADMM-based DMPC approach for controlling the distributions of multi-agent systems under uncertainty.
Furthermore, Bestler and Graichen \cite{bestler2019distributed} and Gobel et al. \cite{gobel2024distributed} adopted ADMM to solve nonlinear DMPC problem that consider the dynamics of coupled subsystems.  
Finally, Tang and Daoutidis \cite{tang2022fast} developed the ELLADA algorithm to address general nonconvex constraints in a distributed setting.  

Although sampling-based optimization methods have been successfully employed in several works to solve nonconvex nonlinear MPC problems, its scalable extension that leverages ADMM for multi-agent MPC has yet to be established.  
Jiang \cite{jiang2024distributed} cast the multi-robot formation navigation task as a probabilistic inference on a graphical model, which is solved using a distributed sampling-based optimizer.  
However, the empirical result in \cite{jiang2024distributed} covers a limited number of agents.   
Furthermore, Dergachev and Yakovlev \cite{dergachev2024model} devised a decentralized MPPI method that constructs a safe distribution by solving a second-order cone programming problem, in which obstacle constraints are not considered.  
Similarly, Streichenberg et al. \cite{streichenberg2023multi} applied decentralized MPPI to the autonomous surface vessel problem.  
However, the method in \cite{streichenberg2023multi} relies on penalty terms and the filtering of unsafe samples, which hinders scalability.  
Wan et al. \cite{wan2021cooperative} offered a distributed stochastic optimal control formulation, but its focus is mainly on cooperative game settings.  

In this paper, we present distributed stochastic search, a scalable sampling-based optimization framework for solving nonconvex nonlinear DMPC problems.  
In Section \ref{sec:preliminary-sampling-based-optimization}, we introduce stochastic search (SS), a theoretical framework that generalizes sampling-based optimization methods.  
In particular, we employ SS to address the nonconvexity of nonlinear MPC problems.  
In Section \ref{sec:distributed-ss}, we formulate a scalable distributed optimization problem and use ADMM to achieve consensus between agents.  
The contributions of this paper are as follows.  
\begin{enumerate}
    \item We present a scalable extension of a generalized sampling-based optimization method that leverages ADMM to solve multi-agent MPC problems involving nonlinear dynamics and nonconvex constraints.  
    \item Through multi-robot navigation simulations, we demonstrate that the proposed method is highly effective in solving nonconvex multi-agent MPC problems compared to a baseline derivative-based method \cite{polik2010interior}.  
\end{enumerate}


\section{\textsc{Sampling-based optimization}}
\label{sec:preliminary-sampling-based-optimization}
\subsection{Variational optimization} 
In variational optimization, we assume that the actual control $\hat{\vu}_k$ is normally distributed as $\hat{\vu}_k \sim \nolinebreak \mathcal{N}(\vu_k, \Sigma)$ \cite{williams2017information}.  
Let $\hat{U} = (\hat{\vu}_{0}, \dots, \hat{\vu}_{T-1})$ be the random sequence of the independent actual inputs and $U = (\vu_{0}, \dots, \vu_{T-1})$ be the command control sequence.  
Let $\Pb_{\hat{U}}$ be the probability measure of the sample space of $\hat{U}$ denoted as $\Omega_{\hat{U}}$ in the uncontrolled system, \textit{i.e.}, $\vu_k = \mathbf{0}$, and $\Qb_{\hat{U}}$ be the measure of open-loop control.  
We denote the probability density functions for $\Pb_{\hat{U}}$ and $\Qb_{\hat{U}}$ by $\vp$ and $\vq$ defined as
\begin{align}
\label{def:marl-sampling-based-opt-pdf}
\vp(\hat{U}) & = \prod\limits_{t=0}^{T-1} Z_{\hat{U}}^{-1} \exp\Big(-\frac{1}{2} \hat{\vu}_t^\top \Sigma^{-1} \hat{\vu}_t \Big), \\
\vq(\hat{U}) & = \prod\limits_{t=0}^{T-1} Z_{\hat{U}}^{-1} \exp\Big(-\frac{1}{2} \big(\hat{\vu}_t \! - \! \vu_t \big)^\top \Sigma^{-1} \big(\hat{\vu}_t \! - \! \vu_t \big)\Big), \\
Z_{\hat{U}} & = \big((2\pi)^{(T)} |\Sigma|\big)^{\frac{1}{2}}.
\end{align}
Let $X_{\hat{U}} = (\vx_{0}, \dots, \vx_{T})$ be the sequence of states following the input sequence $\hat{U}$.  
The state cost function is defined as 
\begin{equation}
J^{\text{state}}(X_{\hat{U}}) = \Phi(\vx_T) + \sum\limits_{t=0}^{T-1} L(\vx_t). 
\end{equation}
Then, the relationship between the free energy of the state cost function and the relative entropy of control distributions can be written as  
\begin{align}
\label{eq:marl-vo-derivation-free-energy}
  \mathcal{F}(J^{\text{state}}(X_{\hat{U}})) &=  -\eta \log\big(\Eb_{\Pb_{\hat{U}}}\big[\exp (-J^{\text{state}}(X_{\hat{U}}) / \eta )\big]\big) \nonumber 
 \\[0.2cm]
 & \leq \Eb_{\Qb_{\hat{U}}}\big[J^{\text{state}}(X_{\hat{U}})\big] + \eta \KL{\Qb_{\hat{U}}}{ \Pb_{\hat{U}}}.
 \end{align}
The right-hand side of the inequality reaches the lower bound, \textit{i.e.}, the free energy by the optimal distribution.  
We describe an optimal distribution $\Qb_{\hat{U}}^*$ in the sample space $\Omega_{\hat{U}}$ that tightens the bound through the optimal density function $\vq^*$ as follows,  
\begin{align}
\label{eq:marl-vo-optimal-control-pdf}
\vq^*(\hat{U}) = Z_{\vq^*}^{-1}\exp(-J^{\text{state}}(X_{\hat{U}}) / \eta)\vp(\hat{U}),
\end{align}
where $Z_{\vq^*}$ is a normalization factor.  
The normalization factor $Z_{\vq^*}$ is calculated using the Monte-Carlo estimate:  
\begin{align}
\label{eq:marl-vo-normalization-factor}
Z_{\vq^*} &= \sum\limits_{m=1}^{M_{\text{sample}}} \exp\bigg(-\frac{1}{\eta} \Big\{J^{\text{state}}\Big(X_{\hat{U}_{\vepsilon}^{(m)}}\Big) \nonumber \\
&- \eta \sum\limits_{t=0}^{T-1} \Big(\frac{1}{2} \vu_t^\top \Sigma^{-1} \big(\vu_t + 2\vepsilon^{(m)}_t \big)\Big) \Big\} \bigg),
\end{align}
where $M_{\text{sample}}$ is the number of samples, $\vepsilon^{(m)}_t$ is the $m$-th noise sample, and $\hat{U}_{\vepsilon}^{(m)} = \big(\vu_{0} + \vepsilon^{(m)}_{0}, \dots, \vu_{T-1} + \vepsilon^{(m)}_{T-1} \big)$.  
Finally, we update the control iteratively as follows.  
\begin{align}
\label{eq:marl-vo-iterative-control-update}
\vu^{(n+1)}_t = \vu^{(n)}_t + \sum\limits_{m=1}^{M_{\text{sample}}} w^{(m)}_{\vepsilon} \vepsilon^{(m)}_t, 
\end{align}
where the superscript $(n)$ represents the update iteration and $w^{(m)}_{\vepsilon}$ is an importance sampling weight.  

\begin{figure*}
   \centering
    \begin{subfigure}[l]{0.485\textwidth}
       \includegraphics[width=\textwidth]{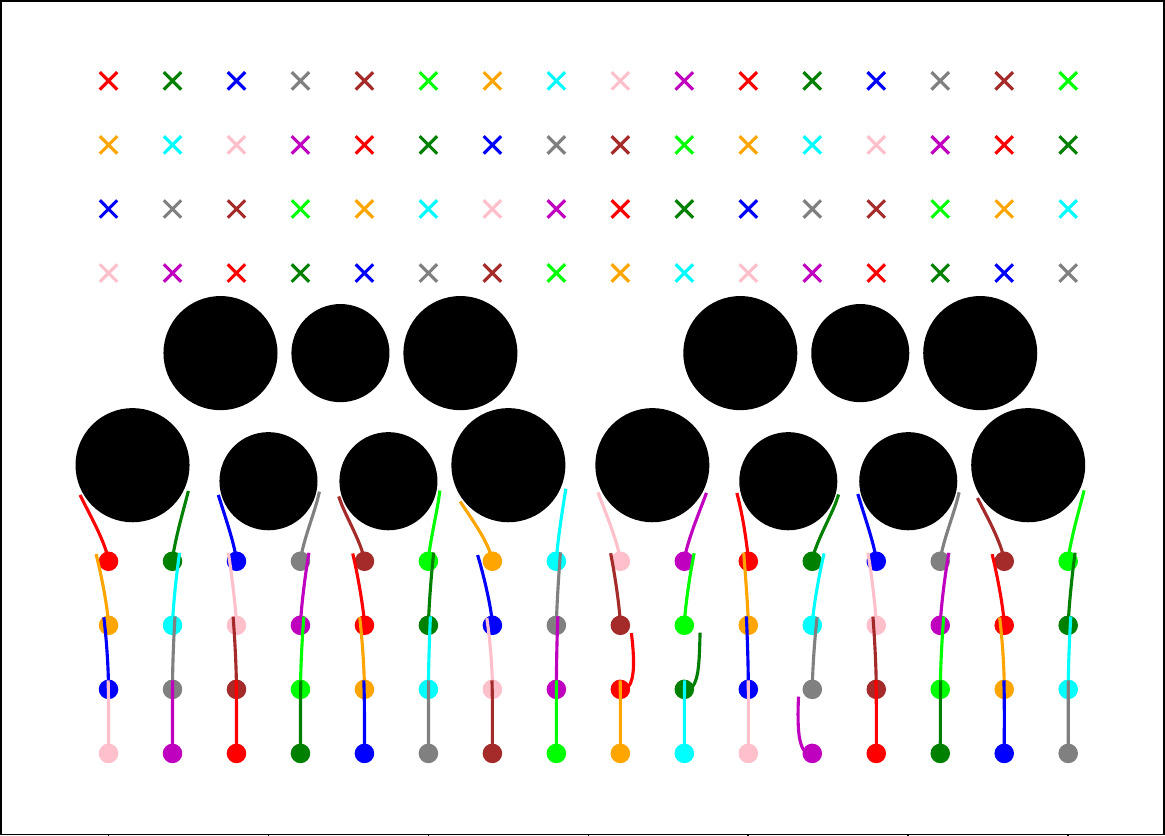}
       \caption{Baseline: Distributed interior point optimizer (IPOPT)}
       \label{fig:dubins-formation-hard1-ipopt}
    \end{subfigure}
    \begin{subfigure}[l]{0.49\textwidth}
       \includegraphics[width=\textwidth]{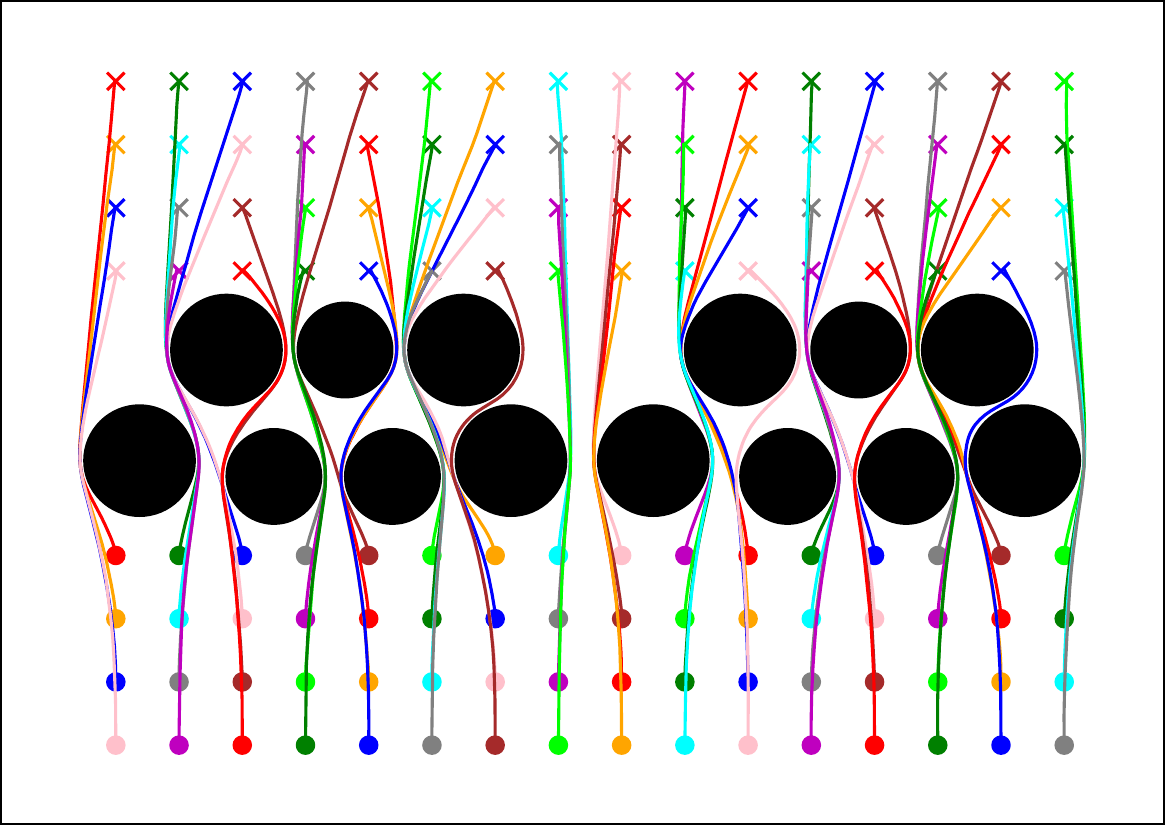}
       \caption{Proposed: Distributed stochastic search (SS)}
       \label{fig:dubins-formation-hard1-ss}
    \end{subfigure}
    \begin{subfigure}[l]{0.488\textwidth}
       \includegraphics[width=\textwidth]{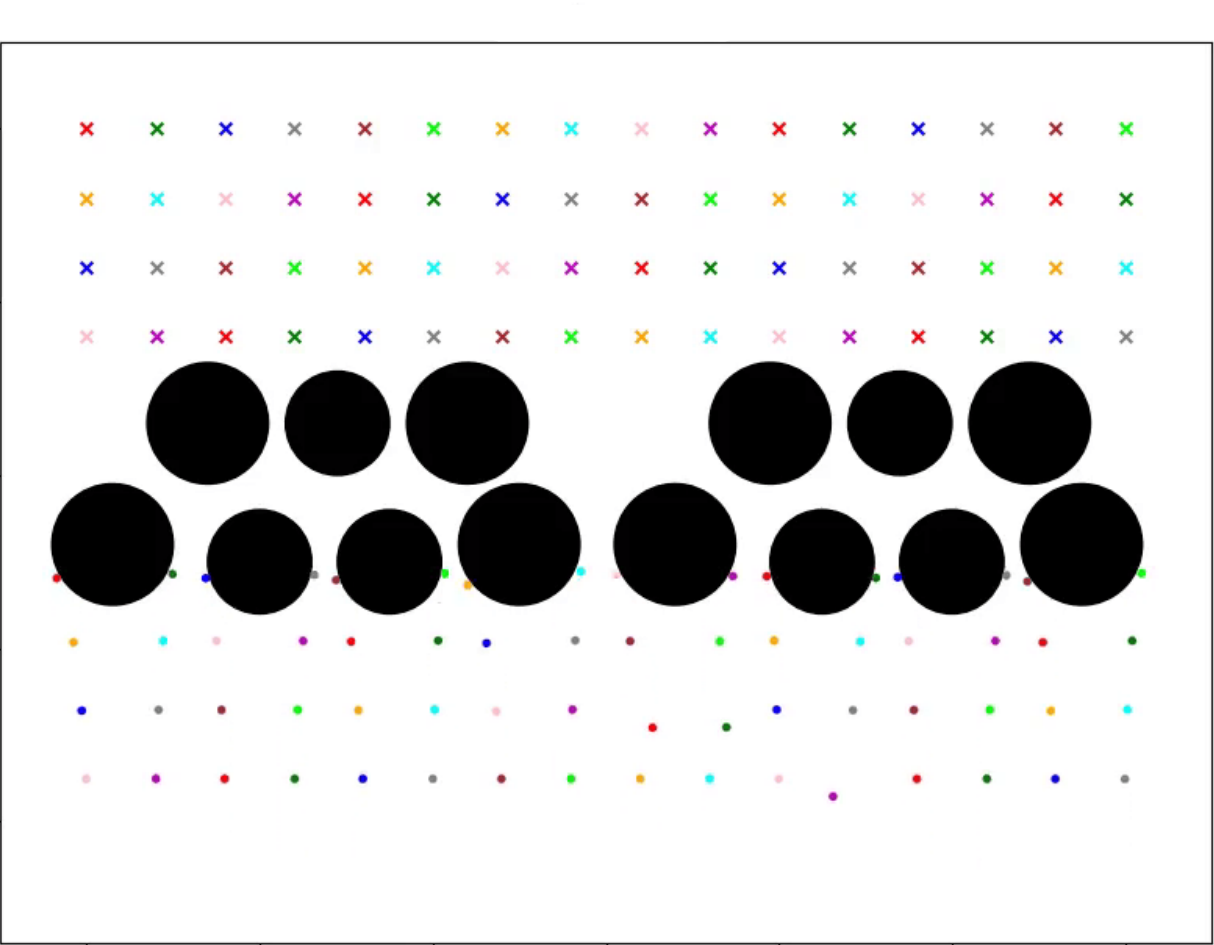}
       \caption{Distributed IPOPT converges to a local minimum}
       \label{fig:dubins-formation-hard1-ipopt-snapshot}
    \end{subfigure}
    \begin{subfigure}[l]{0.488\textwidth}
       \includegraphics[width=\textwidth]{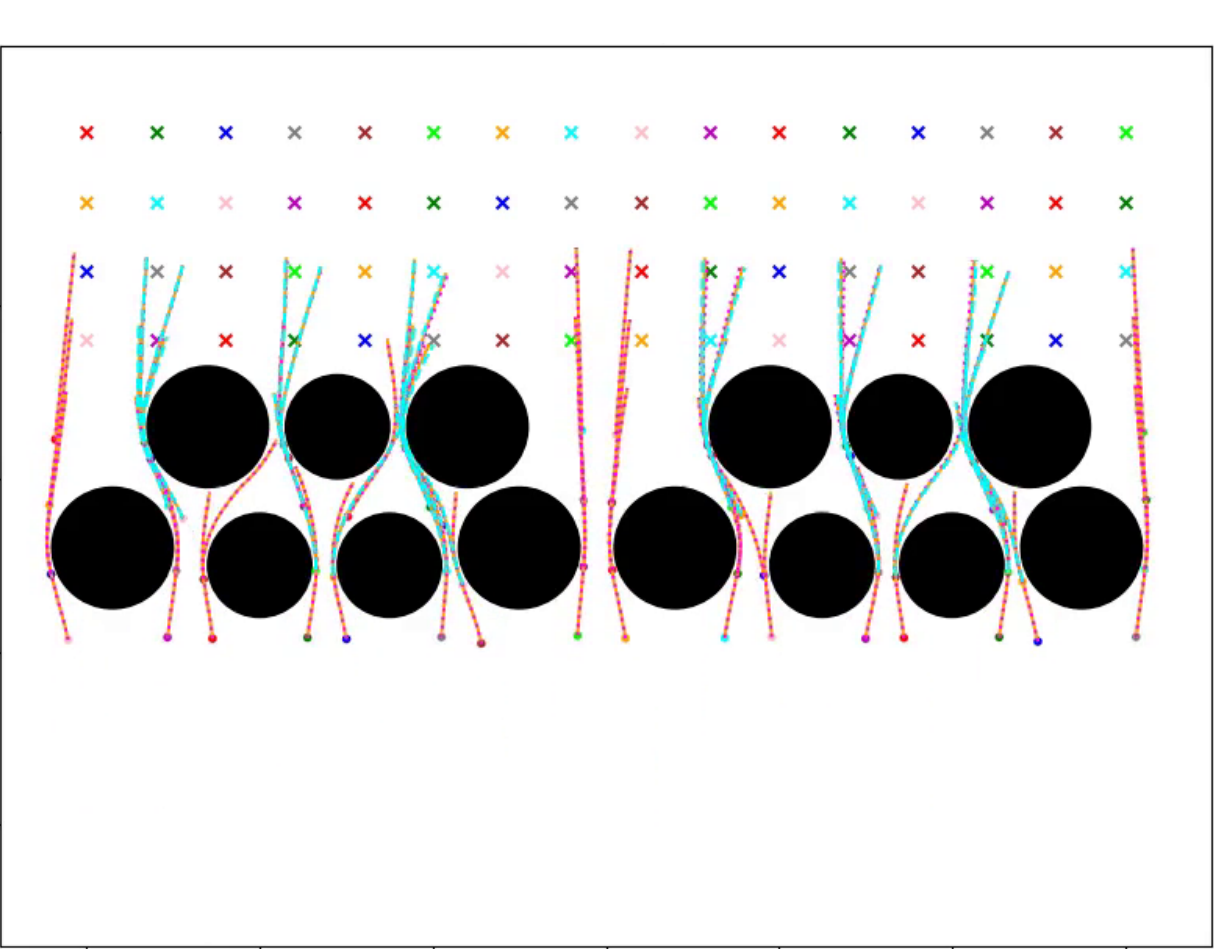}
       \caption{Distributed SS escapes from a local minimum}
       \label{fig:dubins-formation-hard1-ss-snapshot}
    \end{subfigure}
    \caption{64-agent Dubins car formation task visualization: ``$\times$" for target positions, ``$\bullet$" for initial positions, black circles for obstacles, real lines for final trajectories, and dashed lines for the sampled trajectories.  
    (\subref{fig:dubins-formation-hard1-ipopt}), (\subref{fig:dubins-formation-hard1-ipopt-snapshot}): Distributed IPOPT converges to the local minimum.  
    (\subref{fig:dubins-formation-hard1-ss}), (\subref{fig:dubins-formation-hard1-ss-snapshot}): Distributed SS avoids the local minimum by sampling trajectories and completes the task.  
    } 
\label{fig:dubins-formation-hard1}
\end{figure*}

\subsection{Cross-entropy method}
The cross-entropy method (CEM) \cite{rubinstein1997optimization, rubinstein2004cross} considers the estimation of the probability that the sequence of states $X_{\hat{U}}$ results in a state cost $J^{\text{state}}(X_{\hat{U}})$ that does not exceed a given threshold $\gamma$, after following the command control sequence $\hat{U}$ sampled from $\vp(\hat{U})$.  
Formally,   
$$
l^{\text{CEM}} = \Pb_{\hat{U}}\big[J^{\text{state}}(X_{\hat{U}}) \le \gamma\big] = \Eb_{\Pb_{\hat{U}}}\big[\vcalI_{\{J^{\text{state}}(X_{\hat{U}}) \le \gamma\}}\big],
$$
where $\vcalI_{\{\cdot\}}$ is the indicator function.  
One simple way to estimate $l^{\text{CEM}}$ is to use the Monte Carlo estimate: $\frac{1}{M_{\text{sample}}} \sum\limits_{m=1}^{M_{\text{sample}}} \vcalI_{\{J^{\text{state}}(X_{\hat{U}}) \le \gamma\}}$, which is an unbiased estimator of $l^{\text{CEM}}$.  
However, in practice, Monte Carlo simulation requires a large computing resource to accurately estimate $l^{\text{CEM}}$, given that $\{J^{\text{state}}(X_{\hat{U}}) \le \gamma\}$ is a rare event.   
In this context, CEM employs the likelihood ratio estimator based on the importance sampling density as follows:  
\begin{align}
\label{def:marl-cem-importance-sampling}
\hat{l}^{\text{CEM}} = \frac{1}{M_{\text{sample}}} \sum\limits_{m=1}^{M_{\text{sample}}} \vcalI_{\{J^{\text{state}}(X_{\hat{U}}) \le \gamma\}} \frac{\vp(\hat{U})}{\vq(\hat{U})}.  
\end{align}
Hence, the optimal density $\vq^*(\hat{U})$ with zero-variance is 
\begin{align}
\label{eq:marl-cem-optimal-importance-sampling-density}
\vq^*(\hat{U}) = \big(\hat{l}^{\text{CEM}}\big)^{-1} \vcalI_{\{J^{\text{state}}(X_{\hat{U}}) \le \gamma\}} \vp(\hat{U}).
\end{align}
However, $\vq^*(\hat{U})$ depends on the unknown variable $\hat{l}^{\text{CEM}}$.  
In this regard, we substitute $\vq^*(\hat{U})$ with an empirically optimal density $\hat{\vq}^*(\hat{U})$ defined as 
\begin{align}
\label{eq:marl-cem-empirical-optimal-importance-sampling-density}
\hat{\vq}^*(\hat{U}) = Z_{\vq^*}^{-1} \vcalI_{\{J^{\text{state}}(X_{\hat{U}}) \le \gamma\}} \vp(\hat{U}).
\end{align}
In particular, $\hat{\vq}^*(\hat{U})$ is a density for the empirically optimal distribution $\hat{\Qb}_{\hat{U}}^*$.  
The normalization factor $Z_{\vq^*}$ is computed using the Monte Carlo estimate:  
\begin{align}
Z_{\vq^*} = &\sum\limits_{m=1}^{M_{\text{sample}}}\Big[\vcalI_{\{J^{\text{state}}(X^{(n)}_{\hat{U}_{\vepsilon}}) \le \gamma\}} \nonumber \\
& \exp\Big(\sum\limits_{t=0}^{T-1} \Big[ \frac{1}{2} \vu_t^\top \Sigma^{-1} \big(\vu_t + 2\vepsilon^{(m)}_t\big)\Big]\Big)\Big].   
\end{align}
Finally, we update the control iteratively as in (\ref{eq:marl-vo-iterative-control-update}).  


\subsection{Variational inference}
Variational inference (VI) aims to sample from the following posterior density \cite{okada2020variational}:  
\begin{align}
\label{eq:marl-vi-posterior}
\vp(\hat{U}| o=1) = \frac{\vp(o=1 | \hat{U}) \vp(\hat{U})}{\vp(o=1)},
\end{align}
where $o \in \{0, 1\}$, and $o=1$ represents the optimality of a trajectory $X_{\hat{U}}$.  
We describe the optimal probability distribution for the uncontrolled system $\Pb^*_{\hat{U}}$ through the posterior density $\vp(\hat{U}| o=1)$.  
We find the optimal density $\vq^*(\hat{U})$ as 
\begin{align}
\label{eq:marl-vi-optimal-control-density}
&\vq^*(\hat{U}) = \arg\min\limits_{\vq(\hat{U})} \KL{\Qb_{\hat{U}}}{\Pb^*_{\hat{U}}} \nonumber \\
&= \arg\min\limits_{\vq(\hat{U})} \KL{\Qb_{\hat{U}}}{\Pb_{\hat{U}}} - \Eb_{\Qb_{\hat{U}}}\big[\log ( \vp(o=1|\hat{U}) )\big].  
\end{align}
If we choose the log likelihood of VI as $\log ( \vp(o=1|\hat{U}) ) = - \eta^{-1} J^{\text{state}}(X_{\hat{U}})$, then the objective of VI in (\ref{eq:marl-vi-optimal-control-density}) is equivalent to the objective of VO in the inequality in (\ref{eq:marl-vo-derivation-free-energy}).  
In this context, the objective of VI is a generalization of the objective of VO.  
We describe an optimal distribution $\Qb_{\hat{U}}^*$ by the optimal density function $\vq^*({\hat{U}})$ as  
\begin{align}
\label{eq:marl-vi-optimal-control-pdf}
\vq^*(\hat{U}) = Z_{\vq^*}^{-1}\exp(\log \vp(o=1|\hat{U}) )\vp(\hat{U}),
\end{align}
where $Z_{\vq^*}$ is a normalization factor.  

\paragraph{Tsallis variational inference}
The Tsallis variational inference (TVI) \cite{wang2021variational} is a generalization of VI.  
In particular, TVI uses $\mathfrak{q}$-logarithm and $\mathfrak{q}$-exponential to replace the KL divergence with the Tsallis divergence \cite{tsallis1988possible}.  
Hence, the optimal density $\vq^*(\hat{U})$ in TVI is derived as 
\begin{align}
\label{marl-tvi-optimal-control-sequence}
\vq^*(\hat{U}) &= \arg\min\limits_{\vq(\hat{U})} D_{\mathfrak{q}}\big(\Qb_{\hat{U}} \| \Pb^*_{\hat{U}}\big)  \\
&= \arg\min\limits_{\vq(\hat{U})} D_{\mathfrak{q}}\big(\Qb_{\hat{U}} \| \Pb_{\hat{U}}\big) - \Eb_{\Qb_{\hat{U}}}\big[\log_{\mathfrak{q}} ( \vp(o=1|\hat{U}) ) \big]. \nonumber
\end{align}
As in (\ref{eq:marl-vi-optimal-control-pdf}), the optimal density $\vq^*(\hat{U})$ is expressed as 
\begin{align}
\label{eq:marl-tvi-optimal-control-pdf}
\vq^*(\hat{U}) = Z_{\vq^*}^{-1}\exp_{\mathfrak{q}}(\log_{\mathfrak{q}}\vp(o=1|\hat{U}) )\vp(\hat{U}).  
\end{align}

Suppose $\log_{\mathfrak{q}} \vp(o=1|\hat{U}) = - \eta^{-1} J^{\text{state}}(X_{\hat{U}})$.  
It follows that the term $\exp_{\mathfrak{q}}(-\frac{1}{\eta} J^{\text{state}}(X_{\hat{U}}))$ in (\ref{eq:marl-tvi-optimal-control-pdf}) can be reparameterized by setting $\eta = (1 - \mathfrak{q} ) \gamma$:
\begin{align}
\label{eq:marl-tvi-reparameterization}
\exp_{\mathfrak{q}}\big(- J^{\text{state}}(X_{\hat{U}}) / \eta \big) &= \big[1 - (1 - \mathfrak{q}) J^{\text{state}}(X_{\hat{U}}) / \eta\big]^{1/(1 - \mathfrak{q})} \nonumber \\
&= \big[1 - J^{\text{state}}(X_{\hat{U}}) / \gamma \big]^{1/(1 - \mathfrak{q})}.  
\end{align}
We compute (\ref{eq:marl-tvi-reparameterization}) as
\begin{align}
\label{eq:marl-tvi-reparameterization-threshold}
&\big[1 - J^{\text{state}}(X_{\hat{U}}) / \gamma ]^{1/(1 - \mathfrak{q})} &  \\
&= \begin{cases}
\exp\Big(\frac{1}{1 - \mathfrak{q}} \log \Big(1 - J^{\text{state}}(X_{\hat{U}}) / \gamma \Big)\Big), &\text{if } J^{\text{state}} < \gamma \\
0, & \text{otherwise.}
\end{cases}
\nonumber
\end{align}
In (\ref{eq:marl-tvi-reparameterization-threshold}), $\gamma$ is a threshold parameter \cite{wang2021variational}.  
To facilitate the tuning process, we use an \textit{elite fraction} that adjusts $\gamma$ according to the cost scale.  

\subsection{Stochastic search} 
The stochastic search (SS) \cite{boutselis2020constrained, wang2021adaptive, wang2022sampling} finds the optimal density $\vq^*(\hat{U})$ formulated as
\begin{align}
\vq^*(\hat{U}) &= \arg\max\limits_{\vq(\hat{U})} \log\Big(\Eb_{\Pb_{\hat{U}}}\Big[S\big(-J^{\text{state}}(X_{\hat{U}})\big)\Big]\Big) \label{eq:marl-ss-optimal-control-density-1} \\
&= \arg\max\limits_{\vq(\hat{U})} \log\bigg(\Eb_{\Qb_{\hat{U}}}\bigg[\frac{\vp(\hat{U})}{\vq(\hat{U})} S\Big(-J^{\text{state}}(X_{\hat{U}})\Big)\bigg]\bigg).  \label{eq:marl-ss-optimal-control-density-2} 
\end{align}
where $S$ is a monotonically increasing shape function.  
SS generalizes VO, TVI and CEM depending on the shape function $S$: 
i) the exponential function $S(x; \eta) = \exp(\eta^{-1} x)$ in (\ref{eq:marl-ss-optimal-control-density-1}) leads to the VO optimal density function; 
ii) the $\mathfrak{q}$-exponential function $S(x; \eta, \mathfrak{q}) = \exp_{\mathfrak{q}}(\eta^{-1} x)$ with $\mathfrak{q}$-logarithm function corresponds to the optimal density function of TVI in (\ref{marl-tvi-optimal-control-sequence});  
iii) the logistic function $S(x; k, x_0) = \big(1 + \exp(-k(x - x_0))\big)^{-1}$, where $k$ is the logistic growth rate or curve steepness and $x_0$ is the midpoint or $0.5$-quantile of the function closely resembles the optimal density function of CEM in (\ref{eq:marl-cem-empirical-optimal-importance-sampling-density}) with a soft threshold.  
Specifically, $x_0$ corresponds to the threshold $\gamma$ in CEM.  

The SS method offers a general framework for stochastic optimization that updates the parameters of the policy distribution in the exponential family by gradient descent \cite{zhou2014gradient}.  
The update law of SS recovers the update law of VO, TVI, and CEM depending on the shape function $S$.  

We parameterize the control density $\vq(\hat{U})$ by $\zeta$ denoted as $\vq(\hat{U}; \zeta)$.  
In particular, $\zeta$ represents the mean of $\hat{U}$.  
SS updates the parameter $\zeta$ as follows.  
\begin{align}
\label{eq:marl-ss-exponential-family-parameter-update}
\zeta^{(n+1)}_k &= \zeta^{(n)}_k \nonumber \\
&+ \alpha_{\text{SS}}^{(n)} \Bigg(\frac{\Eb\Big[ S(-J^{\text{state}}(X_{\hat{U}}))\Big( \hat{\vu}_k - \zeta^{(n)}_k \Big)\Big]}{\Eb\Big[S(-J^{\text{state}}(X_{\hat{U}}))\Big]}\Bigg),
\end{align}
where $\alpha_{\text{SS}}$ is the step size and the superscript $(n)$ represents the SS iteration.  

\section{\textsc{Distributed stochastic search for multi-agent MPC}}
\label{sec:distributed-ss}
Let $N$ be the number of agents, and $[N] = \{1, \dots, N\}$ denote the index set of $N$.  
Let $\vx$ and $\vu$ be the state and control of the system.  
We describe the discrete-time nonlinear dynamics of the system as 
\begin{align}
\label{eq:marl-mpc-system-dynamics}
\vx_{k + 1} = \vF(\vx_k, \vu_k).  
\end{align}

We define the multi-agent topology as a hypergraph $\calH = ([N], \calE)$, where $[N]$ represents the nodes of agents and $\calE$ is a set of hyperedges.  
A hypergraph is a generalization of a graph that allows any single edge, known as a hyperedge, to connect any number of nodes.  
The neighborhood of agent $i$ at time $k$ is denoted as a hyperedge $\calE_{i,k}$.  
Formally, $\calE_{i,k} = \{j \in [N], j \neq i: d(i, j) \le \delta \}$, where $d$ is a measure that maps the nodes $i$ and $j$ to a real number, and $\delta$ is a threshold.  
In addition, we represent the set of agents that have agent $i$ as their neighbor at time $k$ as $\vec{\calE}_{i,k}$.  

\subsection{Problem formulation}
We assume that each agent $i$ can communicate with another agent $j$ if $i \in \calE_{j,k}$ or $j \in \calE_{i,k}$, and $i \notin \calE_{i,k}$ for all $i \in [N]$ and time steps $k = t_0, \dots, T$.   
The local system dynamics in discrete time is described as 
\begin{align}
\label{def:marl-system-dynamics}
\vx_{i,k+1} = \vF_{i}(\vx_{i,k}, \vu_{i,k}), \forall i \in [N], \forall k.
\end{align}

The total state $\vx$ and control $\vu$ are concatenations of local variables, \textit{i.e.}, $\vx = [\vx_1, \dots, \vx_N]$, $\vu = [\vu_1, \dots, \vu_N]$.  
Similarly, the total cost can be decomposed as 
\begin{equation}
L(\vx_k, \vu_k) = \sum\limits_{i=1}^N L_i(\vx_{i,k}, \vu_{i,k}), ~    
\Phi(\vx_T) = \sum\limits_{i=1}^N \Phi_i(\vx_{i,T}).
\end{equation}
Let $\vX_k^T = (\vx_k, \dots, \vx_T)$ be the sequence of total states and $\vU_k^{T-1} = (\vu_k, \dots, \vu_{T-1})$ be the sequence of total controls.  
Let $\vX_{i,k}^T = (\vx_{i,k}, \dots, \vx_{i,T})$ be the sequence of local states and $\vU_{i, k}^{T-1} = (\vu_{i,k}, \dots, \vu_{i, T-1})$ be the sequence of local controls.  
Let $T_{\text{MPC}}$ be the look-ahead horizon and $T$ be the entire horizon.  
We define the cost function for MPC as 
\begin{align}
\label{def:marl-cost-function}
&J\big(\vX_k^{k+T_{\text{MPC}}}, \vU_k^{k+T_{\text{MPC}} - 1}\big) = \sum\limits_{i=1}^N J_i(\vX_{i,k}^{k+T_{\text{MPC}}}, \vU_{i,k}^{k+T_{\text{MPC}}} - 1)) \nonumber \\
&= \sum\limits_{i=1}^N \Big[ \Phi_i(\vx_{i,k + T_{\text{MPC}}}) +\sum\limits_{t=k}^{k+T_{\text{MPC}} - 1} L_i(\vx_{i,t}, \vu_{i,t}) \Big].
\end{align}

Let $\vg_{i}, \forall i \in [N]$ and $\vh_{i,j}, \forall i,j \in [N], i \neq j$ be functions for individual and inter-agent local constraints, respectively.  

\paragraph{Centralized multi-agent MPC problem}
We formulate the multi-agent MPC problem as follows.  
\begin{problem}[Centralized multi-agent MPC problem]
For each agent $i \in [N]$, find the optimal control sequence $\vU_k^{k+T_{\text{MPC}} - 1}$ such that
\label{prob:centralized-marl}
\begin{align}
\minimize\limits_{\vU_k^{k+T_{\text{MPC}} - 1}} & J(\vX_{k}^{k+T_{\text{MPC}}}, \vU_k^{k+T_{\text{MPC}} - 1})  \nonumber \\
\text{subject to } & \vx_{i,k+1} = \vF_{i}(\vx_{i,k}, \vu_{i,k}), \forall i \in [N], \forall k, \nonumber \\
& \vg_{i}(\vx_{i,k}) \ge 0, \forall i \in [N], \forall k, \nonumber \\
& \vh_{i,j}(\vx_{i,k}, \vx_{j,k}) \ge 0, \forall i, j \in [N], i \neq j, \forall k.  
\end{align}
\end{problem}
Note that Problem \ref{prob:centralized-marl} can be optimized in a centralized manner.  
However, the centralized approach is not scalable as the number of agents increases.  
To improve scalability, we formulate a distributed multi-agent MPC problem.  

\paragraph{Distributed multi-agent MPC problem}
To formulate a distributed approach, we first define augmented variables.  
We define the augmented state and control as $\tilde{\vx}_{i,k} = [\vx_{i,k}, (\vx^i_{j,k})_{j \in \calE_{i,k}}]$ and $\tilde{\vu}_{i,k} = [{\vu}_{i,k}, ({\vu}^i_{j,k})_{j \in \calE_{i,k}}]$ for all $i \in [N]$, where $\vx^i_{j,k}$ is the agent $i$'s copy variable of the agent $j$'s state $\vx_{j,k}$, and $[(\vx^i_{j,k})_{j \in \calE_{i,k}}]$ is a concatenation of vectors $\vx^i_{j,k}, \forall j \in \calE_{i,k}$.  
Let $\tilde{\vF}_i(\tilde{\vx}_{i,k}, \tilde{\vu}_{i,k})$ denote the augmented local dynamics defined as $\tilde{\vF}_i(\tilde{\vx}_{i,k}, \tilde{\vu}_{i,k}) = [\vF_i\big(\vx_{i,k}, \vu_{i,k}), (\vF_j(\vx^i_{j,k}, \vu^i_{j,k}))_{j \in \calE_{i,k}}]$.  
Similarly, $\tilde{\vh}_i(\tilde{\vx}_{i,k})) = [(\vh_{i, j}(\vx_{i,k}, \vx^i_{j,k}))_{j \in \calE_{i,k}}]$.  
For a consensus of local augmented solutions, we denote the global state and control variables as $\vy_k = [\vy_{1,k}, \dots, \vy_{N,k}]$ and $\vz_k = [\vz_{1,k}, \dots, \vz_{N,k}]$.  
To handle constraints, we add the penalty term using the indicator function $\vcalI$ to the cost function.  
Subsequently, we define the augmented cost function:   
\begin{align}
& \tilde{J}_i\big(\vX_{i,k}^{k+T_{\text{MPC}}}, \vU_{i,k}^{k+T_{\text{MPC}} - 1}\big) = J_i\big(\vX_{i,k}^{k+T_{\text{MPC}}}, \vU_{i,k}^{k+T_{\text{MPC}} - 1}\big) \nonumber \\
& + \sum\limits_{t=k}^{k+T_{\text{MPC}} - 1} \Big[\vcalI\Big(\tilde{\vF}_{i}\big(\tilde{\vx}_{i,k}, \tilde{\vu}_{i,k})\Big)\Big] \nonumber \\
&+ \sum\limits_{t=k + 1}^{k+T_{\text{MPC}}} \Big[\vcalI\Big(\vg_{i}(\vx_{i,k})\Big) + \vcalI\Big(\tilde{\vh}_{i}(\tilde{\vx}_{i,k})\Big) \Big].  
\end{align}
To coordinate neighboring agents using ADMM, we reformulate Problem \ref{prob:centralized-marl} using the augmented cost function and the consensus variables as follows.  

\begin{problem}[Distributed multi-agent MPC problem] For each agent $i \in [N]$, find the optimal control sequence $\vU_k^{k+T_{\text{MPC}} - 1}$ such that
\label{prob:marl-augmented-soft-constraints}
\begin{align}
\label{prob:marl-augmented-simplified}
& \minimize\limits_{\vU_{i,k}^{k+T_{\text{MPC}} - 1}} \sum\limits_{i=1}^N \tilde{J}_i\big(\vX_{i,k}^{k+T_{\text{MPC}}}, \vU_{i,k}^{k+T_{\text{MPC}} - 1}\big) \nonumber \\
&\text{subject to } \tilde{\vx}_{i,k} = \tilde{\vy}_{i,k}, \quad \tilde{\vu}_{i,k} = \tilde{\vz}_{i,k}, \forall i \in [N], \forall k,
\end{align}
\end{problem}
where $\tilde{\vy}_{i,k} = [\vy_{i,k}, (\vy_{j,k})_{j \in \calE_{i,k}}]$ and $\tilde{\vz}_{i,k} = [\vz_{i,k}, (\vz_{j,k})_{j \in \calE_{i,k}}]$, $\forall i \in [N], k$.   

\subsection{Consensus optimization using ADMM}
We use ADMM to solve Problem \ref{prob:marl-augmented-soft-constraints} in a distributed manner.  
The augmented Lagrangian is formulated as
\begin{align}
& \calL\big(\tilde{\vx}_k, \tilde{\vu}_k, \tilde{\vy}_k, \tilde{\vz}_k, \boldsymbol{\lambda}_k, \boldsymbol{\xi}_k, \rho, \mu\big) \nonumber \\
& = \sum\limits_{i=1}^N \calL_i\big(\tilde{\vx}_{i,k}, \tilde{\vu}_{i,k}, \tilde{\vy}_{i,k}, \tilde{\vz}_{i,k}, \boldsymbol{\lambda}_{i,k}, \boldsymbol{\xi}_{i,k}, \rho, \mu \big) \nonumber \\
& = \sum_{i=1}^{N} \bigg[ \tilde{J}_i\big(\vX_{i,k}^{k+T_{\text{MPC}}}, \vU_{i,k}^{k+T_{\text{MPC}} - 1}\big) \nonumber \\
&+\sum\limits_{t=k}^{k + T_{\text{MPC}}}  \Big[ \boldsymbol{\lambda}^{\top}_{i,t}  (\tilde{\vx}_{i,t} - \tilde{\vy}_{i,t}) + \frac{\rho}{2}  \| \tilde{\vx}_{i,t} - \tilde{\vy}_{i,t}\|_2^2 \Big] \nonumber \\
&+ \sum\limits_{t=k}^{k + T_{\text{MPC}}-1} \Big[ \boldsymbol{\xi}^{\top}_{i,t}  (\tilde{\vu}_{i,t} - \tilde{\vz}_{i,t}) + \frac{\mu}{2}  \| \tilde{\vu}_{i,t} - \tilde{\vz}_{i,t}\|_2^2 \Big] \bigg].
\end{align}
Let $\calL_i = \calL_i(\tilde{\vx}_{i,k}, \tilde{\vu}_{i,k}, \tilde{\vy}_{i,k}, \tilde{\vz}_{i,k}, \boldsymbol{\lambda}_{i,k}, \boldsymbol{\xi}_{i,k}, \rho, \mu)$ for simplicity.  
The ADMM updates to solve Problem \ref{prob:marl-augmented-soft-constraints} are as follows.  
\paragraph{Local stochastic search update}
The local updates are given by 
\begin{align}
&\big\{\tilde{\vx}_{i,k}^{(l + 1)}, \tilde{\vu}_{i,k}^{(l + 1)}\big\} = \arg\min\limits_{\tilde{\vx}_{i,k}, \tilde{\vu}_{i,k}} \calL_i^{(l)}, \label{eq:marl-global-consensus-admm-variable-update-finite-horizon-mpc}
\end{align}
where the superscript $(l)$ is the ADMM iteration.  
Specifically, SS iteratively updates the augmented mean parameter $\tilde{\zeta}$ from which the augmented control $\tilde{\vu}$ is sampled as
\begin{align}
\label{eq:marl-ss-exponential-family-parameter-update-admm}
\tilde{\zeta}^{(n+1)}_{i,k} &= \tilde{\zeta}^{(n)}_{i,k} \nonumber \\
&+ \alpha_{\text{SS}}^{(n)} \Bigg(\frac{\Eb\Big[ S(-\calL_i^{(l, n)})\Big( \tilde{\vu}_{i,k}^{(l, n)} - \tilde{\zeta}^{(n)}_{i,k} \Big)\Big]}{\Eb\Big[S(-\calL_i^{(l, n)})\Big]}\Bigg), 
\end{align}
for all $i \in [N]$ and $k$, until SS iteration $(n)$ reaches the predefined number of iterations.  

\paragraph{Global and dual updates}
The global variables $\vy, \vz$ are updated as  
\begin{align}
&\vy_{i,k}^{(l + 1)} = \frac{1}{|\vec{\calE}_i|+1} \sum\limits_{j=1}^{|\vec{\calE}_i|+1} \Big[(\vx_{i,k}^j)^{(l+1)} + \frac{ \boldsymbol{\lambda}_{i,k}^{(l)}}{\rho} \Big], \nonumber \\
&\vz_{i,k}^{(l + 1)} = \frac{1}{|\vec{\calE}_i|+1} \sum\limits_{j=1}^{|\vec{\calE}_i|+1} \Big[(\vu_{i,k}^j)^{(l+1)} + \frac{ \boldsymbol{\xi}_{i,k}^{(l)}}{\mu}\Big], \label{eq:marl-global-consensus-admm-variable-update-finite-horizon-global-update} 
\end{align}
where $\rho, \mu > 0$ are the penalty parameters.
Subsequently, the dual variables are updated as
\begin{align}
&\boldsymbol{\lambda}_{i,k}^{(l + 1)} = \boldsymbol{\lambda}_{i,k}^{(l)} + \rho(\tilde{\vx}_{i,k}^{(l + 1)} - \tilde{\vy}_{i,k}^{(l + 1)}), \nonumber \\
&\boldsymbol{\xi}_{i,k}^{(l + 1)} = \boldsymbol{\xi}_{i,k}^{(l)} + \mu(\tilde{\vu}_{i,k}^{(l + 1)} - \tilde{\vz}_{i,k}^{(l + 1)}), \label{eq:marl-global-consensus-admm-variable-update-finite-horizon-dual-update}
\end{align}
for all $i \in [N]$ and $k$.   
In particular, distributed SS iteratively performs the following procedure.  
\begin{enumerate}
    \item \textit{Local SS update:} SS finds the augmented state and control variables that minimizes the augmented Lagrangian $\calL_i^{(l)}$ in (\ref{eq:marl-global-consensus-admm-variable-update-finite-horizon-mpc}).  
    \item \textit{Global update:}  
    The copy variables $\vx_{i,k}^j, \vu_{i,k}^j, \forall j \in \vec{\calE}_i$ are accessed through the local augmented variables $\big\{\tilde{\vx}_{j,k}^{(l + 1)}, \tilde{\vu}_{j,k}^{(l + 1)}\big\}$ collected from each agent $j \in \vec{\calE}_{i, k}$ for all $i \in [N]$.  
    The global variables are then updated as in (\ref{eq:marl-global-consensus-admm-variable-update-finite-horizon-global-update}).   
    \item \textit{Dual update:} The global variables $\vy_{j, k}, \vz_{j, k}$ are collected from each agent $j \in \calE_{i, k}$ for all $i \in [N]$ and the augmented global variables are constructed as $\tilde{\vy}_{i,k}^{(l+1)} = [\vy_{i,k}^{(l+1)}, (\vy_{j,k}^{(l+1)})_{j \in \calE_i}]$ and $\tilde{\vz}_{i,k}^{(l+1)} = [{\vz}_{i,k}^{(l+1)}, ({\vz}_{j, k}^{(l+1)})_{j \in \calE_i}],$ for each agent $i \in [N]$.  
    Subsequently, the dual variables $\boldsymbol{\lambda}_{i,k}^{(l + 1)}$ and $\boldsymbol{\xi}_{i,k}^{(l + 1)}$ are updated as in (\ref{eq:marl-global-consensus-admm-variable-update-finite-horizon-dual-update}).  
\end{enumerate}
Note that the computations in (\ref{eq:marl-global-consensus-admm-variable-update-finite-horizon-mpc})-(\ref{eq:marl-global-consensus-admm-variable-update-finite-horizon-dual-update}) are parallelizable among all agents.

\begin{table*}
\centering
\resizebox{\textwidth}{!}{%
\begin{tabular}{@{}c|ccccccc@{}}
\toprule
\textbf{System task} & \textbf{Type} & \textbf{Optimizer} & \textbf{Task completion rate (\%) $\uparrow$} & \textbf{Inter-agent collisions $\downarrow$} & \textbf{Obstacle collisions $\downarrow$} & \textbf{Trajectory cost $\downarrow$} & \textbf{Computation time (s) $\downarrow$} \\
\midrule
\multirow{8}{*}{Unicycle} &
\multirow{4}{*}{Centralized} 
& MPPI & 3.3 & 0 & 0 & 132612.25 & 467.91 \\
&& CEM & 13.3 & 0 & 0 & 13095.42 & 298.65 \\
&& Tsallis & 13.3 & 0 & 0 & 16045.48 & 289.06 \\
&& IPOPT & local infeasibility &---&---&---&--- \\
\cmidrule(lr){2-8}
& \multirow{4}{*}{Distributed} 
& MPPI & 53.3 & 18 & 0 & 6284.93 & 1574.20 \\
swap && CEM & 33.3 & 0 & 0 & 6240.87 & 1629.57 \\
30 agents && \textbf{Tsallis} & \textbf{53.3} & \textbf{0} & \textbf{0} & 6268.47 & 1644.43 \\
&& IPOPT & local infeasibility & --- & --- & --- & --- \\
\cmidrule(lr){1-8}

\multirow{8}{*}{Unicycle} & 
\multirow{4}{*}{Centralized} 
& MPPI & 2.08 & 0 & 2 & 725162.76 & 738.06 \\
&& CEM & 8.33 & 0 & 0 & 84349.20 & 934.22 \\
&& Tsallis & 8.33 & 0 & 0 & 82795.00 & 727.32 \\
&& IPOPT & local infeasibility & --- & --- & --- & --- \\
\cmidrule(lr){2-8}
& \multirow{4}{*}{Distributed} 
& MPPI & 87.5 & 0 & 0 & 25390.31 & 2803.63 \\
formation && \textbf{CEM} & \textbf{89.6} & \textbf{0} & \textbf{0} & 25474.02 & 2811.86 \\
48 agents && \textbf{Tsallis} & \textbf{89.6} & \textbf{0} & \textbf{0} & 25531.90 & 2873.56 \\
&& IPOPT & restoration failed & --- & --- & --- & --- \\
\cmidrule(lr){1-8}

\multirow{8}{*}{Dubins} & 
\multirow{4}{*}{Centralized} 
& MPPI & 0 & 0 & 23 & 21712401.53 & 2628.54 \\
&& CEM & 0 & 0 & 9 & 3139040.42 & 3079.45 \\
&& Tsallis & 0 & 0 & 0 & 1270955.87 & 3017.25 \\
&& IPOPT & time limit & --- & --- & --- & --- \\
\cmidrule(lr){2-8}
& \multirow{4}{*}{Distributed} 
& MPPI & 95.3 & 0 & 0 & 473763.80 & 3547.87 \\
formation && \textbf{CEM} & \textbf{100} & \textbf{0} & \textbf{0} & 441553.26 & 3894.83 \\
64 agents && \textbf{Tsallis} & \textbf{100} & \textbf{0} & \textbf{0} & 442100.65 & 3624.89 \\
&& IPOPT & local infeasibility & --- & --- & --- & --- \\
\cmidrule(lr){1-8}

\multirow{8}{*}{Quadcopter} & 
\multirow{4}{*}{Centralized} 
& MPPI & 0 & 0 & 4 &  --- &  627.08 \\
&& CEM & 0 & 0 & 0 & 38643578219.54 & 528.91 \\
&& Tsallis & 4.2 & 0 & 0 & 3780913792.34 & 625.79 \\
&& IPOPT & time limit & --- & --- & --- & --- \\
\cmidrule(lr){2-8}
& \multirow{4}{*}{Distributed} 
& MPPI & 12.5 & 130 & 0 & 65616485.14 & 15450.91 \\
formation && CEM & 25 & 32 & 0 & 38677991.44 & 12676.57 \\
24 agents && \textbf{Tsallis} & \textbf{25} & \textbf{6} & \textbf{0} & 37443003.73 & 12303.18 \\
&& IPOPT & local infeasibility & --- & --- & --- & --- \\
\bottomrule
\end{tabular}%
}
\caption{Simulations for comparing stochastic search (SS) and interior point optimizer (IPOPT) in centralized and distributed optimization on unicycle, Dubins, and quadcopter systems.  
Model predictive path integral (MPPI) \cite{williams2017information}, cross entropy method (CEM) \cite{rubinstein2004cross}, and Tsallis MPPI \cite{wang2021variational} are generalized by SS.
}  
\label{tab:hard1-comparison}
\end{table*}

\begin{figure}
   \centering
       \includegraphics[width=0.42\textwidth]{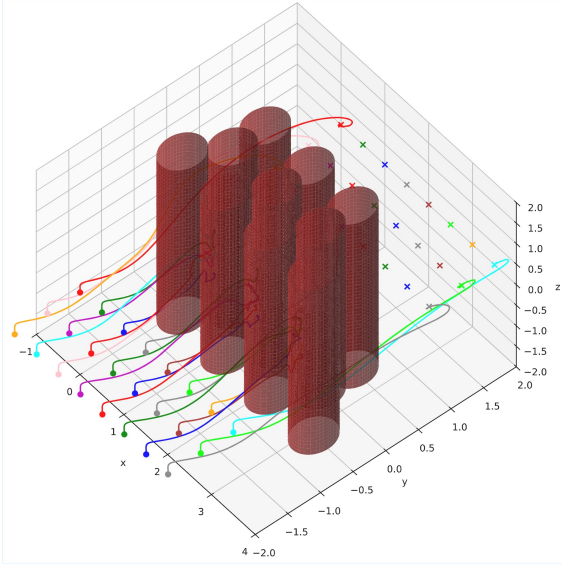}
       \label{fig:dubins-formation-hard2-quadcopter}
    \caption{Distributed stochastic search (proposed method) trajectories on 24-agent quadcopter formation task: ``$\times$" for target positions, ``$\bullet$" for initial positions, red purple columns for obstacles, and real lines for trajectories.
    } 
\label{fig:quadcopter-formation-hard2}
\end{figure}

\section{\textsc{Simulation results}}
We compare the performance of our method (SS) with the interior point optimizer (IPOPT) \cite{wachter2006implementation}, which is a successful convex optimization NLP solver.  
In particular, we evaluate methods on multi-robot navigation tasks, in which each agent is governed by nonlinear dynamics.  
For distributed IPOPT, the multi-robot navigation tasks are implemented using the code developed in \cite{mallick2024multi}.  

Table \ref{tab:hard1-comparison} shows that distributed optimization methods outperform centralized optimization methods in all cases.  
In all systems and tasks, distributed Tsallis MPPI consistently outperforms the other methods, while IPOPT consistently fails due to local infeasibility or high computational complexity.  
The local infeasibility output of IPOPT refers to the status that the solver has converged to a local minimizer for constraint violations in the L1 norm but is infeasible for the problem \cite{wachter2006implementation}.  

In Fig. \ref{fig:dubins-formation-hard1}, distributed IPOPT fails due to local infeasibility during optimization, while distributed SS successfully reaches the goal in the 64-agent multi-car navigation task.  
In Fig. \ref{fig:quadcopter-formation-hard2}, we test the proposed method in highly nonconvex problems, in which we intentionally block the path of agents using obstacles in the quadcopter formation task.   
Despite the highly nonconvex setup, distributed Tsallis MPPI shows 25\% task completion rate with the lowest number of collisions.  

\subsection{Multi-robot navigation}
In multi-robot navigation, we use the following quadratic cost function defined as 
\begin{align}
\label{eq:marl-exp-cost-function}
&J_i(\vX_{i, k}^{k+T_{\text{MPC}}}, \vU_{i, k}^{k+T_{\text{MPC}} - 1}) \\
&= (\vx_{i, k+T_{\text{MPC}}} - \vx_{i}^{\text{ref}})\T \vQ^f_i (\vx_{i, k+T_{\text{MPC}}} - \vx_{i}^{\text{ref}}) \nonumber \\
&+ \sum\limits_{t=k}^{k+T_{\text{MPC}}}\big[(\vx_{i, t} - \vx_{i}^{\text{ref}})\T \vQ_i (\vx_{i, t} - \vx_{i}^{\text{ref}}) + \vu_{i, t}\T \vR_i \vu_{i, t}\big].
\nonumber
\end{align}
The same trajectory cost coefficients are used for both SS and IPOPT.    
For SS, we add the penalty $10^{10}$ as an indicator function to the cost function in (\ref{eq:marl-exp-cost-function}) for each collision.  
For IPOPT, we impose hard constraints to avoid obstacle/inter-agent collisions.  
Additionally, we impose jump prevention constraints in IPOPT, \textit{i.e.}, state constraints to prevent agents from jumping over other agents or obstacles in a single time interval.  
To avoid non-differentiability of the Euclidean distance at zero, we apply the offset of $10^{-8}$ to the L2 norm in IPOPT constraints.  

The entire time horizon is set $T = 280$ for the formation task and $T = 180$ for the swap task with the time interval $dt = 0.02$, and the look-ahead horizon is set $T_{\text{MPC}} = 30$.  
For SS, we use $M_{\text{sample}} = 1024$ samples.  

\paragraph{Unicycle}
In a two-dimensional unicycle system, the state vector is $\vx = [x^{\text{pos}}, y^{\text{pos}}, \theta]^\top \in \Rb^3$, and the control vector is $\vu = [v, \omega]^\top \in \Rb^2$.  
To be specific, $x^{\text{pos}}, y^{\text{pos}}$ are $xy$-coordinates, $\theta$ is an orientation, $v$ is the linear velocity, and $\omega$ is the angular velocity of the unicycle.  
The unicycle system is governed by the following dynamics:
\begin{align}
\label{eq:marl-exp-unicycle-dynamics}
\mathbf{F}(\mathbf{x}_k, \mathbf{u}_k) = [v_{k} \cos(\theta_k\big), v_k \sin(\theta_k), \omega_k]^\top.
\end{align}
The following cost coefficients are used to optimize the unicycle system: $\vQ_i = \operatorname{diag}([0.8, 0.8, 0]), \vQ^f_i = \operatorname{diag}([10, 10, 0])$, and $\vR_i = \operatorname{diag}([0.1, 0.05])$.  
For SS, the standard deviation of the sampling distribution is assigned as $\Sigma = \operatorname{diag}([0.8, 0.8])$.  
For IPOPT, the convergence tolerance is set $10^{-3}$.  
The number of centralized [Tsallis, CEM, MPPI] iterations for the i) formation task is [50, 50, 50], and ii) swap task is [50, 50, 80].  
The number of [ADMM, SS] iterations for distributed SS is [50, 2].  
The number of [ADMM, consensus] iterations for distributed IPOPT is [10, 20].  
The ADMM penalty coefficient for distributed SS is 0.2 and distributed IPOPT is 0.5.  

\paragraph{Dubins}
The state vector of the Dubins car system is defined as $\vx = [x^{\text{pos}}, y^{\text{pos}}, \theta, v]^\top \in \Rb^4$ and the control vector is defined as $\vu = [\omega, a]^\top \in \Rb^2$, where $a$ is the linear acceleration.  
The Dubins car system is governed by the following dynamics:
\begin{equation}
\label{eq:marl-exp-dubins-dynamics}
\mathbf{F}(\mathbf{x}_k, \mathbf{u}_k) = [v_{k} \cos(\theta_k), v_k \sin(\theta_k), \omega_k, a_k ]^\top.
\end{equation}
The following cost coefficients are used to optimize the Dubins system: $\vQ_i = \vQ^f_i = \operatorname{diag}([6, 6, 0, 1])$ and $\vR_i = \operatorname{diag}([0.1, 0.05])$.  
For SS, the standard deviation of the sampling distribution is set  $\Sigma = \operatorname{diag}([0.8, 0.8])$.  
For IPOPT, the convergence tolerance is set $10^{-3}$.  
The number of centralized [Tsallis, CEM, MPPI] iterations for the formation task is [70, 70, 60].  
The number of [ADMM, SS] iterations for distributed SS is [20, 2].  
The number of [ADMM, consensus] iterations for distributed IPOPT is [10, 20].  

\paragraph{Quadcopter}
In a quadcopter system, the state vector is defined as $\vx = [x^{\text{pos}}, y^{\text{pos}}, z^{\text{pos}}, v^x, v^y, v^z, \phi^{\text{roll}}, \theta^{\text{pitch}}, \psi^{\text{yaw}}, v^{\phi}, v^{\theta}, v^{\psi}]^\top \in \Rb^{12}$ 
and the control vector is defined as $\vu = [f^{\text{thrust}}, a^{\phi}, a^{\theta}, a^{\psi}]^\top \in \Rb^4$.  
To be specific, $v^x, v^y, v^z$ are linear velocities, $\phi^{\text{roll}}, \theta^{\text{pitch}}, \psi^{\text{yaw}}$ are Euler angles, $v^{\phi}, v^{\theta}, v^{\psi}$ are Euler angle angular velocities, $f^{\text{thrust}}$ is the thrust, and $a^{\phi}, a^{\theta}, a^{\psi}$ are Euler angle angular accelerations.  
We use the quadcopter dynamics derived in \cite{sabatino2015quadrotor}.  
The following cost coefficients are used to optimize the quadcopter system: $ \vQ_i = \operatorname{diag}([1000, 1000, 1000, 100, 100, 100, 1000, 1000, 1000, \allowbreak 100, 100, 100])$, $\vQ^f_i = 100 \cdot \vQ_i$, and $\vR_i = \operatorname{diag}([1, 1, 1, 1])$.  
For SS, the standard deviation of the sampling distribution is assigned as $\Sigma = \operatorname{diag}([2, 1, 1, 1])$.  
For IPOPT, the convergence tolerance is set $10^{-8}$.  
The number of centralized [Tsallis, CEM, MPPI] iterations for the formation task is [60, 50, 60].  
The number of ADMM iterations for distributed [Tsallis, CEM, MPPI] is [40, 40, 50].  
The number of SS iterations for distributed SS is 2.  
The number of [ADMM, consensus] iterations for distributed IPOPT is [10, 20].  
The ADMM penalty coefficient for distributed [Tsallis, CEM, MPPI] is [0.85, 1, 1].    
The ADMM penalty coefficient for distributed IPOPT is 0.5.    

\section{\textsc{Conclusion}}
In this paper, we have proposed distributed stochastic search, a generalized and scalable sampling-based optimization framework for multi-agent model predictive control.  
The proposed method utilizes sampling-based optimization to address nonconvexity and ADMM to achieve scalability over an increasing number of agents.  
In multi-robot navigation simulations, the proposed method shows a scalable performance to overcome nonconvexity of the problem compared to the distributed baseline method that uses interior point optimizer, a derivative-based nonlinear programming solver.  

In future research, improving the constraint handling capability of our method would be an interesting direction, since tight constraint satisfaction is not yet theoretically grounded due to its stochastic nature.  
Another promising direction would be to consider uncertainty through chance-constrained \cite{farina2016stochastic} or robust optimization formulations \cite{bertsimas2021probabilistic, abdul2025scalable}.  
Finally, it would be interesting to explore distributed learning-to-optimize approaches \cite{saravanos2025deep} for accelerating the convergence of the proposed optimizer. 


\addtolength{\textheight}{-0cm}   

\bibliographystyle{IEEEtran}
\bibliography{IEEEabrv,root}

\end{document}

%% file: commands.tex
\newcommand{\tr}{\rm{tr}}
\newcommand{\Var}{Var}
\newcommand{\Cov}{Cov}
 
\newcommand{\SaS}{{\mathcal{S}\alpha\mathcal{S}}}
\newcommand{\xb}{\mathbf{x}}
\newcommand{\ub}{\mathbf{u}}
\newcommand{\wb}{\mathbf{w}}
\newcommand{\rd}{{\mathrm d}}
\newcommand{\rT}{{\mathrm T}}
\newcommand{\calU}{{\cal{U}}}
\newcommand{\calV}{{\cal{V}}}
\newcommand{\calO}{{\cal{O}}}
\newcommand{\calS}{{\mathcal{S}}}
\newcommand{\calI}{{\mathcal{I}}}
\newcommand{\calE}{{\mathcal{E}}}
\newcommand{\calA}{{\mathcal{A}}}
\newcommand{\calB}{{\mathcal{B}}}
\newcommand{\calG}{{\mathcal{G}}}
\newcommand{\T}{^\top}
\newcommand{\Sr}{\mathcal{S}_r}
 
\newcommand{\rpartial}{{\mathrm{\partial}}}
\newcommand{\rinf}{{\mathrm \inf}}
\newcommand{\rp}{{\mathrm p}}
\newcommand{\rdelta}{{\mathrm \delta}}
\newcommand{\rtr}{{\mathrm{tr}}}
\newcommand{\rP}{{\mathrm{P}}}
\newcommand{\rvec}{{\mathrm{vec}}}
\newcommand{\rtau}{{\mathrm{\tau}}}
\newcommand{\va}{{\bf a}}
\newcommand{\vb}{{\bf b}}
\newcommand{\vc}{{\bf c}}
\newcommand{\vx}{{\bf x}}
\newcommand{\vy}{{\bf y}}
\newcommand{\vz}{{\bf z}}
\newcommand{\vp}{{\bf p}}
\newcommand{\dvx}{{\bf dx}}
\newcommand{\vdu}{{\bf du}}
\newcommand{\vdy}{{\bf dy}}
\newcommand{\vq}{{\bf q}}
\newcommand{\vr}{{\bf r}}
\newcommand{\vs}{{\bf  s}}
\newcommand{\vf}{{\bf f}}
\newcommand{\vg}{{\bf g}}
\newcommand{\vu}{{\bf u}}
\newcommand{\vv}{{\bf  v}}
\newcommand{\vh}{{\bf h}}
\newcommand{\vl}{{\bf l}}
\newcommand{\vdx}{{\bf dx}}
\newcommand{\vdw}{{\bf dw}}
\newcommand{\vw}{{\bf w}}
\newcommand{\bm}{{\bf 1}}
 
\newcommand{\vA}{{\bf A}}
\newcommand{\vB}{{\bf B}}
\newcommand{\vC}{{\bf C}}
\newcommand{\vD}{{\bf D}}
\newcommand{\vE}{{\bf E}}
\newcommand{\vF}{{\bf F}}
\newcommand{\vG}{{\bf G}}
\newcommand{\vH}{{\bf H}}
\newcommand{\vI}{{\bf I}}
\newcommand{\vJ}{{\bf J}}
\newcommand{\vK}{{\bf K}}
\newcommand{\vL}{{\bf L}}
\newcommand{\vM}{{\bf M}}
\newcommand{\vN}{{\bf N}}
\newcommand{\vO}{{\bf O}}
\newcommand{\vP}{{\bf P}}
\newcommand{\vQ}{{\bf Q}}
\newcommand{\vR}{{\bf R}}
\newcommand{\vS}{{\bf S}}
\newcommand{\vT}{{\bf T}}
\newcommand{\vU}{{\bf U}}
\newcommand{\vV}{{\bf V}}
\newcommand{\vW}{{\bf W}}
\newcommand{\vX}{{\bf X}}
\newcommand{\vY}{{\bf Y}}
 
\newcommand{\vOmega}{{\bf \Omega}}
\newcommand{\vxi}{{\mbox{\boldmath$\xi$}}}
\newcommand{\vpi}{{\mbox{\boldmath$\pi$}}}
\newcommand{\vdomega}{{\bf d\omega}}
\newcommand{\vlambda}{{\mbox{\boldmath$\lambda$}}}
\newcommand{\vBamma}{{\mbox{\boldmath$\Gamma$}}}
\newcommand{\VTheta}{{\mbox{\boldmath$\Theta$}}}
\newcommand{\VPhi}{{\mbox{\boldmath$\Phi$}}}
\newcommand{\vphi}{{\mbox{\boldmath$\phi$}}}
\newcommand{\VPsi}{{\mbox{\boldmath$\Psi$}}}
\newcommand{\vepsilon}{{\mbox{\boldmath$\epsilon$}}}
\newcommand{\vSigma}{{\mbox{\boldmath$\Sigma$}}}
\newcommand{\valpha}{{\mbox{\boldmath$\alpha$}}}
\newcommand{\vmu}{{\mbox{\boldmath$\mu$}}}
\newcommand{\vbeta}{{\mbox{\boldmath$\beta$}}}
\newcommand{\vomega}{{\mbox{\boldmath$\omega$}}}
\newcommand{\vtau}{{\mbox{\boldmath$\tau$}}}
\newcommand{\vdtau}{{\mbox{\boldmath$d\tau$}}}
\newcommand{\vtheta}{{\mbox{\boldmath$\theta$}}}\newcommand{\dataset}{{\cal D}}
\newcommand{\fracpartial}[2]{\frac{\partial #1}{\partial  #2}}
\newcommand{\vcalS}{{\mbox{\boldmath$\cal{S}$}}}
\newcommand{\vcalU}{{\mbox{\boldmath$\cal{U}$}}}
\newcommand{\vcalD}{{\mbox{\boldmath$\cal{D}$}}}
\newcommand{\vcalJ}{{\mbox{\boldmath$\cal{J}$}}}
\newcommand{\vcalE}{{\mbox{\boldmath$\cal{E}$}}}
\newcommand{\vcalF}{{\mbox{\boldmath$\cal{F}$}}}
\newcommand{\calH}{{\cal H}}
\newcommand{\calX}{{\cal X}}
\newcommand{\calF}{{\cal F}}
\newcommand{\calL}{{\cal L}}
\newcommand{\calP}{{\cal P}}
\newcommand{\calN}{{\cal N}}
\newcommand{\calM}{{\cal M}}
\newcommand{\vcalL}{{\mbox{\boldmath$\cal{L}$}}}
\newcommand{\vcalZ}{{\mbox{\boldmath$\cal{Z}$}}}
\newcommand{\vcalG}{{\mbox{\boldmath$\cal{G}$}}}
\newcommand{\vcalN}{{\mbox{\boldmath$\cal{N}$}}}
\newcommand{\vcalM}{{\mbox{\boldmath$\cal{M}$}}}
\newcommand{\vcalH}{{\mbox{\boldmath$\cal{H}$}}}
\newcommand{\vcalC}{{\mbox{\boldmath$\cal{C}$}}}
\newcommand{\vcalO}{{\mbox{\boldmath$\cal{O}$}}}
\newcommand{\vcalP}{{\mbox{\boldmath$\cal{P}$}}}
\newcommand{\vcalB}{{\mbox{\boldmath$\cal{B}$}}}
\newcommand{\vcalA}{{\mbox{\boldmath$\cal{A}$}}}
\newcommand{\vcalg}{{\mbox{\boldmath$\cal{g}$}}}
\newcommand{\vcalI}{{\mbox{\boldmath$\cal{I}$}}}
\newcommand{\argmax}{\operatornamewithlimits{argmax}}
\newcommand{\argmin}{\operatornamewithlimits{argmin}}
 
\newcommand{\vzeta}{{\bf \zeta}}
\newcommand{\veta}{{\bf \eta}}
\newcommand{\vgamma}{{\bf \gamma}}
\newcommand{\vGamma}{{\bf \Gamma}}
\newcommand{\vDelta}{{\bf \Delta}}
 
\newcommand{\Qb}{\mathbb{Q}}
\newcommand{\Pb}{\mathbb{P}}
\newcommand{\Xspace}{\mathcal{X}}
\newcommand{\snew}{\mathbf{x}^\prime}
\newcommand{\s}{\mathbf{x}}
\newcommand{\Pas}[2]{\mathcal{P}\left({#1}|{#2}\right)}
\newcommand{\Pol}[2]{\mathcal{U}\left({#1}|{#2}\right)}
\newcommand{\PolOpt}[2]{\mathcal{U}^{*}\left({#1}|{#2}\right)}
\newcommand{\PolPasT}[1]{\mathcal{P}\left({#1}\right)}
\newcommand{\PolT}[1]{\mathcal{U}\left({#1}\right)}
\newcommand{\traj}{\mathbf{X}}
\newcommand{\KL}[2]{\mathcal{D}_{KL}\left({#1}\parallel {#2}\right)}
\newcommand{\TKL}[3]{\mathcal{D}_{TKL_{#1}}\left({#2}\parallel {#3}\right)}
\newcommand{\f}[2]{\mathcal{D}_f\left({#1}\parallel {#2}\right)}
\newcommand{\costt}[1]{\mathcal{J}\left({#1}\right)}
\newcommand{\logb}[1]{\log\left({#1}\right)}
\newcommand{\expb}[1]{\exp\left({#1}\right)}
\newcommand{\val}[2]{\mathit{V}_{#2}\left({#1}\right)}
\newcommand{\zt}[2]{\Phi_{#2}\left({#1}\right)}
\newcommand{\nZ}[2]{\mathcal{G}_t[\Phi]\left({#1}\right)}
 
\newcommand{\Nb}{\mathbb{N}}
\newcommand{\Rb}{\mathbb{R}}
\newcommand{\Sb}{\mathbb{S}}
\newcommand{\Eb}{\mathbb{E}}
\newcommand{\Ib}{\mathbb{I}}
\newcommand{\CVaR}{\text{CVaR}}
\newcommand{\VaR}{\text{VaR}}
\newcommand{\ExP}[2]{\Eb_{{#1}}{\left[#2\right]}}
 
\newcommand{\rf}{{\mathrm f}}
\newcommand*\samethanks[1][\value{footnote}]{\footnotemark[#1]}
 
\newcommand{\dom}{\operatorname{dom}}
 
\newcommand*\dif{\mathop{} \mathrm{d}}
 
\newcommand*\da{\dif a}
\newcommand*\db{\dif b}
\newcommand*\dc{\dif c}
\newcommand*\dd{\dif d}
\newcommand*\de{\dif e}
\newcommand*\df{\dif f}
\newcommand*\dg{\dif g}
\newcommand*\di{\dif i}
\newcommand*\dk{\dif k}
\newcommand*\dl{\dif l}
\newcommand*\dm{\dif m}
\newcommand*\dn{\dif n}
\newcommand*\dq{\dif q}
\newcommand*\dr{\dif r}
\newcommand*\ds{\dif s}
\newcommand*\dt{\dif t}
\newcommand*\du{\dif u}
\newcommand*\dv{\dif v}
\newcommand*\dw{\dif w}
\newcommand*\dx{\dif x}
\newcommand*\dy{\dif y}
\newcommand*\dz{\dif z}
 
\newcommand*\dA{\dif A}
\newcommand*\dB{\dif B}
\newcommand*\dC{\dif C}
\newcommand*\dD{\dif D}
\newcommand*\dE{\dif E}
\newcommand*\dF{\dif F}
\newcommand*\dG{\dif G}
\newcommand*\dH{\dif H}
\newcommand*\dI{\dif I}
\newcommand*\dJ{\dif J}
\newcommand*\dK{\dif K}
\newcommand*\dL{\dif L}
\newcommand*\dM{\dif M}
\newcommand*\dN{\dif N}
\newcommand*\dO{\dif O}
\newcommand*\dP{\dif P}
\newcommand*\dQ{\dif Q}
\newcommand*\dR{\dif R}
\newcommand*\dS{\dif S}
\newcommand*\dT{\dif T}
\newcommand*\dU{\dif U}
\newcommand*\dV{\dif V}
\newcommand*\dW{\dif W}
\newcommand*\dX{\dif X}
\newcommand*\dY{\dif Y}
\newcommand*\dZ{\dif Z}

\newcommand{\minimize}{\mathop{\operatorname{minimize~}}}
 
\newcommand{\todo}[1]{\textcolor{red}{\textbf{TODO:} #1}}
\newcommand{\note}[1]{\textcolor{blue}{\textbf{NOTE:} #1}}